\documentclass{article}%
\usepackage{amsmath}
\usepackage{amsfonts}
\usepackage{amssymb}
\usepackage{graphicx}%
\usepackage{mathrsfs}
\DeclareMathOperator*{\esssup}{ess\,sup}
\providecommand{\U}[1]{\protect \rule{.1in}{.1in}}
\newtheorem{theorem}{Theorem}[section]

\newtheorem{lemma}[theorem]{Lemma}

\newtheorem{remark}[theorem]{Remark}

\newenvironment{proof}[1][Proof]{\noindent \textbf{#1.} }{\  \rule{0.5em}{0.5em}}
\begin{document}

\title{Multidimensional quadratic BSDEs with separated generators}
\author{Asgar Jamneshan \thanks{Department of Mathematics and Statistics, University of Konstanz} \and{Michael Kupper\thanks{Department of Mathematics and Statistics, University of Konstanz}}\and{Peng Luo \thanks{Department of Mathematics, ETH Zurich, 8092 Zurich, Switzerland. }}}
\date{September 11, 2017}

\maketitle

\begin{abstract}
We consider multidimensional quadratic BSDEs with bounded and unbounded terminal conditions. We provide sufficient conditions which guarantee existence and uniqueness of solutions. In particular, these conditions are satisfied if the terminal condition or
the dependence in the system are small enough.
\end{abstract}


\textbf{MSC-classification}: 60H10, 60H30.
\section{Introduction}
Backward stochastic differential equations (BSDEs) are introduced in Bismut \cite{B}.
A BSDE is an equation of the form
\begin{equation}\label{bsde}
Y_t=\xi+\int_{t}^{T}g(s,Y_s,Z_s)ds-\int_{t}^{T}Z_sdW_s, \quad t\in[0,T],
\end{equation}
where $W$ is a $d$-dimensional Brownian motion,  the terminal condition $\xi$ is an $n$-dimensional random variable, and $g:\Omega\times [0,T]\times \mathbb{R}^n \times \mathbb{R}^{n\times d}\to \mathbb{R}^n$ is the generator. A solution consists of a pair of predictable processes $(Y,Z)$ with values in $\mathbb{R}^n$
and $\mathbb{R}^{n\times d}$, called the value and control process, respectively. The first existence and uniqueness result for BSDEs with an $L^2$-terminal condition and a generator satisfying a Lipschitz growth condition is due to  Pardoux and Peng~\cite{PP}. In case that the generator satisfies a quadratic growth condition in the control $z$, the situation is more involved and a general existence theory  does not exist. Frei and dos Reis \cite{FR} and Frei \cite{F} provide counterexamples which show that multidimensional quadratic BSDEs may fail to have a global solution. In the one-dimensional case the existence of quadratic BSDE
is shown by Kobylanski \cite{Ko} for bounded terminal conditions, and by Briand and Hu \cite{BH,BH1} for unbounded terminal conditions.
Solvability results for superquadratic BSDEs are discussed in Delbaen et al.~\cite{DHB}, see also Masiero and Richou \cite{MR}, Richou \cite{R} and Cheridito and Nam  \cite{CN1}.

The focus of the present work lies on multidimensional quadratic BSDEs.
In case that the terminal condition is small enough the
existence and uniqueness of a solution was first shown by Tevzadze \cite{T}. Cheridito and Nam \cite{CN} and Hu and Tang \cite{HT} obtain local solvability on $[T-\varepsilon,T]$ for some $\varepsilon>0$ of systems of BSDEs with subquadratic generators and diagonally quadratic generators respectively, which under additional
assumptions on the generator can be extended to global solutions.
Cheridito and Nam \cite{CN} provide solvability for Markovian quadratic BSDEs and projectable quadratic BSDEs . Frei \cite{F} introduced the notion of split solution and studied the existence of multidimensional quadratic BSDEs by considering a special kind of terminal condition. In Bahlali et al.~\cite{BEH} existence is shown when the generator $g(s,y,z)$ is strictly subquadratic in $z$ and satisfies some monotonicity condition.

Our results are motivated by the recent work of Hu and Tang \cite{HT}. We focus on the solvability of multidimensional quadratic BSDEs with generators that are independent of the value process. More precisely, we study the coupled system of quadratic BSDEs
\begin{equation}\label{S2}
Y^{i}_{t}=\xi^i+\int_{t}^{T}g^i(s,Z_s)ds-\int_{t}^{T}Z^{i}_sdW_s, \quad t\in[0,T], \; i=1,\dots,n
\end{equation}
where $g^i(s,z)=f^{i}(s,z^i)+h^i(s,z)$. As suggested by Hu and Tang \cite{HT}, any solution of \eqref{S2} is a fixed point of the mapping $I$ defined by $z\mapsto I(z):=Z$ where $Z$ is given by the decoupled system of BSDEs
\begin{equation}\label{S3}
Y^{i}_{t}=\xi^i+\int_{t}^{T}f^i(s,Z^i_s)+h^i(s,z_s)ds-\int_{t}^{T}Z^{i}_sdW_s, \quad t\in[0,T], \; i=1,\dots,n.
\end{equation}
In order to solve the decoupled BSDEs \eqref{S3}, we adapt the techniques of Briand and Hu \cite{BH} and provide existence and uniqueness
for the one-dimensional BSDEs with unbounded terminal conditions. Under the quadratic growth assumptions $|f^i(t,z^i)-f^i(t,\bar{z}^i)|\leq \theta_i(|z^i|+|\bar{z}^i|)|z^i-\bar{z}^i|$ and $|h^i(t,z)-h^i(t,\bar{z})|\leq\vartheta_i(|z|+|\bar{z}|)|z-\bar{z}|$
for some constants $\theta_i\ge 0$ and $\vartheta_i\ge 0$, $i=1,\dots,n$, we state in Theorem \ref{thm2} the sufficient conditions
\eqref{thmeq1}--\eqref{thmeq2} on the interplay
between the size of $\theta_i$, $\vartheta_i$ and $\xi^i$, which guarantees a unique BMO-solution of \eqref{S2}.
Similar results are given in Theorem \ref{thm3} for unbounded terminal conditions. We realize that the function
\begin{equation*}
u(\gamma,x)=\frac{1}{(2\gamma)^2}\left(e^{2\gamma x}-1-2\gamma x\right)
\end{equation*}
introduced in Briand and Hu \cite{BH} plays an essential role in this setting. Indeed, the function $u$ is a Lyapunov function which is used in the recent work of Xing and \v{Z}itkovi\'{c} \cite{XZ} where the global solvability is obtained for a large class of multidimensional quadratic BSDEs in the Markovian setting. The function $u$ allows us to consider different cases:
(i) if $\theta_i=0$,
the conditions \eqref{thmeq1}--\eqref{thmeq2} are satisfied under smallness of the terminal condition, in which case we recover the condition in Tevzadze \cite{T} for generators that are independent of the value process, or the condition in  Kramkov and Pulido \cite{KP}; (ii) for any $\theta_i$ and $\xi^i$ the the conditions \eqref{thmeq1}--\eqref{thmeq2} are always satisfied if $\vartheta_i$ are small enough. %We give some explanation in Remark \ref{rm}.
Our results could be applied to market making problems (see  Kramkov and Pulido \cite{KP}), nonzero-sum risk-sensitive stochastic differential games (see El Karoui and Hamad\`{e}ne \cite{EH}, Hu and Tang \cite{HT}) and non-zero sum differential games of BSDEs (see Hu and Tang \cite{HT1}).

The paper is organized as follows. In Section \ref{M}, we state the setting and main results. An auxiliary result for 1-dimensional quadratic BSDEs is presented in the Appendix \ref{A}.

\section{Preliminaries and main results}\label{M}
Let $W=(W_t)_{t\geq 0}$ be a $d$-dimensional Brownian motion on a probability space $(\Omega, {\cal F}, P)$. Let $(\mathcal{F}_t)_{t\geq 0}$ be the augmented filtration generated by $W$. Throughout, we fix a $T\in (0,\infty)$. We endow $\Omega \times [0,T]$ with the predictable $\sigma$-algebra $\mathcal{P}$ and $\mathbb{R}^n$ with its Borel $\sigma$-algebra $\mathcal{B}(\mathbb{R}^n)$. Equalities and inequalities between random variables and processes are understood in the $P$-a.s. and $P\otimes dt$-a.e. sense, respectively. For two real numbers $a,b\geq 0$, their minimum is denoted by $a\wedge b$. The Euclidean norm is denoted by $|\cdot|$ and $\|\cdot\|_\infty$ denotes the $L^\infty$-norm. Let $\mathcal{S}^\infty(\mathbb{R}^n)$ be the space of all $n$-dimensional continuous adapted processes such that
\begin{equation*}
\|Y\|_{\infty}:=\| \sup_{0\leq t\leq T} |Y_t| \|_{\infty} < \infty.
\end{equation*}
Let $\mathcal{T}$ be the set of all stopping times with values in $[0,T]$. For any uniformly integrable martingale $M$ with $M_0=0$, define
\begin{equation*}
\|M\|_{BMO}:=\sup_{\tau\in\mathcal{T}}\|E[|\langle M\rangle_T-\langle M\rangle_\tau||\mathcal{F}_\tau]^{\frac{1}{2}}\|_\infty.
\end{equation*}
The class $\{M:~\|M\|_{BMO}<\infty\}$ is denoted by $BMO$, which is written as $BMO(P)$ when it is necessary to indicate the underlying probability measure $P$. For $(\alpha\cdot W)_t:=\int_0^t\alpha_sdW_s$ in $BMO$, the corresponding stochastic exponential is denoted by $\mathcal{E}_{t}(\alpha\cdot W)$. We recall a classical result on $BMO$ spaces (see \cite[Theorem 3.6]{Ka}).
\begin{lemma}\label{PP}
Let $a\cdot W\in BMO$ be such that $\|a\cdot W\|_{BMO}\leq K$ for some $K\geq 0$, and $\tilde{P}$ be given by $\frac{d\tilde{P}}{dP}:=\mathcal{E}_T(a\cdot W)$, under which $\tilde{W}=W-\int_0^{\cdot}a_sds$ is a Brownian motion. Then for every $b\cdot W\in BMO$, there exist two constants $l(K)$ and $L(K)$ only depending on $K$ such that
\begin{equation*}
 l(K)\|b\cdot W\|^2_{BMO} \leq \| b\cdot \tilde{W}\|^2_{BMO(\tilde{P})} \leq L(K) \| b\cdot W\|^2_{BMO}.
 \end{equation*}
\end{lemma}
In the following we define
\begin{equation*}
h(K):=\frac{L^2(K)}{l(K)}
\end{equation*}
for $K\geq 0$. By \cite[Theorem 2]{CM}, we can choose $l(K)=\frac{2}{(K+\sqrt{2})^2}$. When $K<\sqrt{2}$, it follows from \cite[Lemma 3.1]{KXZ} that $L(K)=\frac{2}{(\sqrt{2}-K)^2}$.
As in Briand and Hu \cite{BH}, for every $\gamma\ge 0$, we consider the function $u(\gamma,\cdot):\mathbb{R}_+\to\mathbb{R}_+$ defined by
\begin{equation*}
u(\gamma,x)=\frac{1}{(2\gamma)^2}\left(e^{2\gamma x}-1-2\gamma x\right).
\end{equation*}
Notice that $u(0,x)=\frac{x^2}{2}$. It is straightforward to check that $x \mapsto u(\gamma,|x|)$ is $\mathcal{C}^2$ and $u''(\gamma,x)-2\gamma u'(\gamma,x)=1$. We consider the  multidimensional quadratic BSDE \eqref{S3}
where $f^i:\Omega\times [0,T]\times \mathbb{R}^{d} \to \mathbb{R}$ and $h^i:\Omega\times [0,T]\times \mathbb{R}^{n\times d} \to \mathbb{R}$ are $\mathcal{P}\otimes \mathcal{B}(\mathbb{R}^{n\times d})$-measurable.
We stipulate the following conditions. For each $i=1,\dots,n$, there are constants $\theta_i\geq 0$ and $\vartheta_i\geq 0$ such that
\begin{itemize}
\item[(A1)] $\xi^i\in L^\infty(\mathcal{F}_T)$;
\item[(A2)] $f^i(t,0)=0$ and $|f^i(t,z^i)-f^i(t,\bar{z}^i)|\leq \theta_i(|z^i|+|\bar{z}^i|)|z^i-\bar{z}^i|$;
\item[(A3)] $h^i(t,0)=0$ and $|h^i(t,z)-h^i(t,\bar{z})|\leq\vartheta_i(|z|+|\bar{z}|)|z-\bar{z}|$.
\end{itemize}
\begin{remark}
 The subsequent results can be extended to the case where $f^i(t,0)$ and $h^i(t,0)$ are bounded, as well as
\begin{equation*}
|f^i(t,z^i)-f^i(t,\bar{z}^i)|\leq \theta_i(1+|z^i|+|\bar{z}^i|)|z^i-\bar{z}^i|,\quad |h^i(t,z)-h^i(t,\bar{z})|\leq\vartheta_i(1+|z|+|\bar{z}|)|z-\bar{z}|.
\end{equation*}
The proof is along the same line of reasoning, however, for the sake of readability we assume (A2) and (A3).
\end{remark}
In the following we state our main existence and uniqueness result.
\begin{theorem}\label{thm2} Suppose (A1)-(A3) hold and the conditions
  \begin{align}\label{thmeq1}
  &\bigg\|\sum_{i=1}^{n}u(\theta_i,|\xi^i|)\bigg\|_\infty+c\bigg\|\sum_{i=1}^{n}\vartheta_iu'\left(\theta_i,\frac{1}{2\theta_i}\log\left(\frac{E[e^{2\theta_i|\xi^i|}|\mathcal{F}_t]}{1-2\theta_i\vartheta_ic}\right)\right)\bigg\|_{\infty}\leq\frac{1}{2}c
  \end{align}
    \begin{align}\label{thmeq2}
 &\sum_{i=1}^{n}16\vartheta_i^2ch(2\theta_i\sqrt{c})<1.
  \end{align}
  are satisfied, where
  \begin{align*}
  c:=\left(\min_{1\leq i\leq n}\frac{1}{4\theta_i\vartheta_i}\right)\wedge4\bigg\|\sum_{i=1}^{n}u(\theta_i,|\xi^i|)\bigg\|_\infty.
  \end{align*}
Then the system \eqref{S2} admits a unique solution $(Y,Z)$ with $Y\in\mathcal{S}^{\infty}(\mathbb{R}^n)$ and $\|Z\cdot W\|^2_{BMO}\leq c$.
\end{theorem}
\begin{proof}
Let $z\cdot W\in BMO$ be such that $\|z\cdot W\|^2_{BMO}\leq c$, by Lemma \ref{le}, the BSDE
\begin{equation*}
Y^{i}_{t}=\xi^i+\int_{t}^{T}f^i(s,Z^{i}_{s})+h^i(s,z_s)ds-\int_{t}^{T}Z^{i}_sdW_s,
\end{equation*}
admits a unique solution $(Y^i,Z^i\cdot W)\in\mathcal{S}^{\infty}(\mathbb{R})\times BMO$ with
\begin{align*}
|Y^i_t|&\leq \|\xi^i\|_\infty+\frac{1}{2\theta_i}\log\left(\frac{1}{1-2\theta_i\vartheta_i\|z\cdot W\|^2_{BMO}}\right)\\
&\leq \|\xi^i\|_\infty+\frac{1}{2\theta_i}\log\left(\frac{1}{1-2\theta_i\vartheta_ic}\right)
\end{align*}
and
\begin{align*}
\frac{1}{2}E\left[\int_t^T|Z^i_s|^2ds\bigg|\mathcal{F}_t\right]&\leq E[u(\theta_i,|\xi^i|)|\mathcal{F}_t]+E\left[\int_t^T\vartheta_iu'\left(\theta_i,\frac{1}{2\theta_i}\log\left(\frac {E[e^{2\theta_i|\xi^i|}|\mathcal{F}_s]}{1-2\theta_i\vartheta_ic}\right)\right)z^2_sds\bigg|\mathcal{F}_t\right].
\end{align*}
Hence,
\begin{align*}
\frac{1}{2}\|Z\cdot W\|^2_{BMO}
&\leq \bigg\|\sum_{i=1}^nu(\theta_i,|\xi^i|)\bigg\|_{\infty}+\bigg\|\sum_{i=1}^{n}\vartheta_iu'\left(\theta_i,\frac{1}{2\theta_i}\log\left(\frac{E[e^{2\theta_i|\xi^i|}|\mathcal{F}_t]}{1-2\theta_i\vartheta_ic}\right)\right)\bigg\|_{\infty}c.
\end{align*}
Consider the candidate set
\begin{equation*}
M=\left\{z: \|z\cdot W\|^2_{BMO}\leq c\right\}.
\end{equation*}
For $z\in M$, define $I(z)=Z$, where $Z$ is given by
\begin{align*}
&Y^i_t=\xi^i+\int_t^Tf^i(s,Z^i_s)+h^i(s,z_s)ds-\int_t^TZ^i_sdW_s, ~~i=1,\dots,n.
\end{align*}
 By \eqref{thmeq1}, $I$ maps $M$ to itself. For $z,\bar{z}\in M$, let $Z=I(z)$ and $\bar{Z}=I(\bar{z})$. Denote $\delta Z:=Z-\bar{Z}$, $\delta z:=z-\bar{z}$, $\delta Y:=Y-\bar{Y}$. Then one has
\begin{align*}
\displaystyle \delta Y^i_t&=\int_t^Tf^i(s,Z^i_s)-f^i(s,\bar{Z}^i_s)+h^i(s,z_s)-h^i(s,\bar{z}_s)ds-\int_t^T\delta Z^i_sdW_s\\
\displaystyle &=\int_t^Ta^i_s\delta Z^i_s+h^i(s,z_s)-h^i(s,\bar{z}_s)ds-\int_t^T\delta Z^i_sdW_s\\
\displaystyle&=\int_t^Th^i(s,z_s)-h^i(s,\bar{z}_s)ds-\int_t^T\delta Z^i_sd\tilde{W}^i_s,
\end{align*}
where $\tilde{W}^i_t:=W_t-\int_0^t a^i_sds$ is a Brownian motion under the equivalent probability measure $\frac{d\tilde{P}^i}{dP}=\mathcal{E}_T(a^i\cdot W)$, and the process $a$ satisfies $|a^i_s|\leq \theta_i|Z^i_s+\bar{Z}^i_s|$. Using It\^{o}'s formula to $|\delta Y^i_t|^2$, taking conditional expectation with respect to $\mathcal{F}_t$ and $\tilde{P}^i$ and using $2ab\leq \frac{1}{2}a^2+2b^2$, one has
\begin{align*}
|\delta Y^i_t|^2+\tilde{E}^i\left[\int_t^T|\delta Z^i_s|^2ds\bigg|\mathcal{F}_t\right]\leq \frac{1}{2}\|\delta Y^i\|^2_{\infty}+2|\vartheta_i|^2\left(\tilde{E}^i\left[\int_t^T\left(|z_s|+|\bar{z}_s|\right)|\delta z_s|ds\bigg|\mathcal{F}_t\right]\right)^2.
\end{align*}
Noting that
\begin{align*}
\frac{1}{2}\left(\|\delta Y^i\|_{\infty}^{2}+\|\delta Z^{i}\cdot W\|_{BMO(\tilde{P}^i)}^{2}\right)
&\leq\|\delta Y^i\|_{\infty}^{2} \vee \|\delta Z^i\cdot W\|_{BMO(\tilde{P}^i)}^{2} \\
&\leq\esssup_{(\omega,t)}\left\{|\delta Y^i_{t}|^2+\tilde{E}^i\left[\int_{t}^{T}|\delta Z^i_s|^2ds\big|\mathcal{F}_t\right]\right\},
\end{align*}
it follows from H\"{o}lder's inequality that
\begin{align*}
\|\delta Z^{i}\cdot W\|_{BMO(\tilde{P}^i)}^{2}&\leq 4|\vartheta_i|^2\esssup_{(\omega,t)}\left(\tilde{E}^i\left[\int_t^T\left(|z_s|+|\bar{z}_s|\right)^2\bigg|\mathcal{F}_t\right]\tilde{E}^i\left[\int_t^T|\delta z_s|^2ds\bigg|\mathcal{F}_t\right]\right),
\end{align*}
and therefore
\begin{align*}
\|\delta Z^i\cdot \tilde{W}\|^2_{BMO(\tilde{P}^i)}\leq 8|\vartheta_i|^2\left(\|z_s\cdot\tilde{W}\|^2_{BMO(\tilde{P}^i)}+\|\bar{z}_s\cdot\tilde{W}\|^2_{BMO(\tilde{P}^i)}\right)\|\delta z_s\cdot\tilde{W}\|^2_{BMO(\tilde{P}^i)}.
\end{align*}
By Lemma \ref{PP}, one has
\begin{align*}
l(2|\theta_i|\sqrt{c})\|\delta Z^i\cdot W\|^2_{BMO}&\leq 8L^2(2\theta_i\sqrt{c})\vartheta_i^2\left(\|z_s\cdot W\|^2_{BMO}+\|\bar{z}_s\cdot W\|^2_{BMO}\right)\|\delta z_s\cdot W\|^2_{BMO}\\
&\leq 16L^2(2\theta_i\sqrt{c})\vartheta_i^2c\|\delta z_s\cdot W\|^2_{BMO}.
\end{align*}
Therefore, it holds that
\begin{align*}
\|\delta Z^i\cdot W\|^2_{BMO}\leq 16h(2\theta_i\sqrt{c})\vartheta_i^2c\|\delta z_s\cdot W\|^2_{BMO}.
\end{align*}
By \eqref{thmeq2}, $I$ is a contraction mapping and the statement follows from the Banach fixed point theorem.
\end{proof}
\begin{remark}
In case $\theta_i=0$ for all $i=1,\dots,n$, it follows that $c=2\|\sum\limits_{i=1}^n|\xi^i|^2\|_{\infty}$. Therefore, the conditions \eqref{thmeq1}--\eqref{thmeq2} reduce to
  \begin{align*}
  &\bigg\|E\left[\sum_{i=1}^{n}\vartheta_i|\xi^i||\mathcal{F}_t\right]\bigg\|_{\infty}\leq\frac{1}{4},\\
  &32\left(\sum_{i=1}^{n}\vartheta_i^2\right)\|\sum_{i=1}^n|\xi^i|^2\|_{\infty}<1.
  \end{align*}
  Noting that
  \begin{align*}
  \bigg\|E\left[\sum_{i=1}^{n}\vartheta_i|\xi^i|\bigg|\mathcal{F}_t\right]\bigg\|_{\infty}&\leq  \bigg\|E\left[\sqrt{\sum_{i=1}^{n}\vartheta_i^2}\sqrt{\sum_{i=1}^n|\xi^i|^2}\bigg|\mathcal{F}_t\right]\bigg\|_{\infty}\\
  &\leq\sqrt{\sum_{i=1}^{n}\vartheta_i^2}\sqrt{\bigg\|\sum\limits_{i=1}^n|\xi^i|^2\bigg\|_{\infty}},
  \end{align*}
  it is enough to assume $32\left(\sum_{i=1}^{n}\vartheta_i^2\right)\|\sum_{i=1}^n|\xi^i|^2\|_{\infty}<1$ to guarantee the solvability of BSDE \eqref{S2}, improving the assumption in \cite[Theorem 11.6]{NT}, see also Tevzadze \cite{T}.
\end{remark}
\begin{remark}\label{R1}
Given $\xi^i$ and $\theta_i$ for $i=1,\dots,n$, the conditions \eqref{thmeq1}--\eqref{thmeq2} are satisfied for sufficiently small $\vartheta_i$, $i=1,\dots,n.$ Indeed, for $\vartheta_i$ small enough, one has $c=4\|\sum_{i=1}^{n}u(\theta_i,|\xi^i|)\|_\infty$, so that the conditions \eqref{thmeq1}--\eqref{thmeq2} reduce to
\begin{align*}
  &\bigg\|\sum_{i=1}^{n}\vartheta_iu'\left(\theta_i,\frac{1}{2\theta_i}\log\left(\frac{E[e^{2\theta_i|\xi^i|}|\mathcal{F}_t]}{1-2\theta_i\vartheta_ic}\right)\right)\bigg\|_{\infty}\leq\frac{1}{4},\\
 &\sum_{i=1}^{n}\left(16\vartheta_i^2ch(2\theta_i\sqrt{c})\right)<1.
  \end{align*}
\end{remark}
Multidimensional quadratic BSDEs with unbounded terminal conditions are first studied in Frei \cite{F} and recently in Kramkov and Pulido \cite{KP}. By using the result in the Appendix \ref{A}, Theorem \ref{thm2} extends to the unbounded case.
\begin{theorem}\label{thm3}
Assume that $\xi^i=E[\xi^i]+\int_0^Tv^i_sdW_s$, where $v^i\cdot W\in BMO$, for all $i=1,\dots,n$.
Denote
\begin{equation*}
C=\left(\min_{1\leq i\leq n}\frac{1}{8\theta_i\vartheta_i}\right)\wedge\sum_{i=1}^n\frac{4\|v^i\cdot W\|^2_{BMO}}{1-16\theta_i\|v^i\|_{BMO}}.
\end{equation*}
Suppose (A2)-(A3) hold and the conditions
  \begin{align}\label{thmeq3}
  &\|v^i\cdot W\|_{BMO}<\frac{1}{16\theta_i},
  \end{align}
  \begin{align}\label{thmeq4}
  &\sum_{i=1}^n\left(\frac{1}{1-16\theta_i\|v^i\|_{BMO}}+\frac{1}{1-4\theta_i\vartheta_iC}\right)\|v^i\cdot W\|^2_{BMO}\leq\frac{1}{2}C,
  \end{align}
  \begin{align}\label{thmeq5}
  &\sum_{i=1}^n\left(16\vartheta_i^2Ch(2\theta_i\sqrt{C})\right)<1,
  \end{align}
  are satisfied. Then the system \eqref{S2} admits a unique solution $(Y,Z)$ with $\|Z\cdot W\|^2_{BMO}\leq C$.
\end{theorem}
\begin{proof}
For $z\cdot W\in BMO$ satisfying $\|z\cdot W\|^2_{BMO}\leq C$, from Lemma \ref{le} it follows that the BSDE
\begin{equation*}
Y^{i}_{t}=\xi^i+\int_{t}^{T}f^i(s,Z^{i}_{s})+h^i(s,z_s)ds-\int_{t}^{T}Z^{i}_sdW_s,
\end{equation*}
admits a unique solution $(Y^i,Z^i\cdot W)$ such that
\begin{align*}
\frac{1}{2}\|Z^i\cdot W\|^2_{BMO}&\leq\left(\frac{1}{1-16\theta_i\|v^i\cdot W\|_{BMO}}+\frac{1}{1-4\theta_i\|z\cdot W\|^2_{BMO}}\right)\|v^i\cdot W\|^2_{BMO}.
\end{align*}
Let
\begin{equation*}
M=\left\{z: \|z\cdot W\|^2_{BMO}\leq C\right\}.
\end{equation*}
For $z\in M$, define $I(z)=Z$, where $Z$ is given by
\begin{align*}
Y^i_t=\xi^i+\int_t^Tf^i(s,Z^i_s)+h^i(s,z_s)ds-\int_t^TZ^i_sdW_s,~~~~ i=1,\dots,n.
\end{align*}
 By similar arguments as in Theorem \ref{thm2}, under the conditions \eqref{thmeq3}-\eqref{thmeq5}, it follows that $I$ is a contraction mapping.
\end{proof}
\begin{remark}
In case $\theta_i=0$ for all $i=1,\dots,n$, it follows that $C=\sum\limits_{i=1}^n4\|v^i\cdot W\|^2_{BMO}$. Thereby the conditions \eqref{thmeq3}-\eqref{thmeq5} reduce to
    \begin{align*}
  &\sum_{i=1}^n2\|v^i\cdot W\|^2_{BMO}\leq\frac{1}{2}C,\\
  &\sum_{i=1}^n\left(16\vartheta_i^2\right)C<1.
  \end{align*}
  Therefore we have to assume $64\left(\sum_{i=1}^n\vartheta_i^2\right)\left(\sum_{i=1}^n\|v^i\cdot W\|^2_{BMO}\right)<1$ to guarantee the solvability of BSDE \eqref{S2}, providing the same bound as in Kramkov and Pulido \cite{KP}.
\end{remark}
\begin{remark}
In line with Remark \ref{R1}, given $\xi^i$ and $\theta_i$ satisfying \eqref{thmeq3}, the conditions \eqref{thmeq4}--\eqref{thmeq5} are satisfied for sufficiently small $\vartheta_i$, $i=1,\dots,n.$
\end{remark}

%\begin{remark}\label{rm}
%Since we can allow $\theta_i=0$, one can put $f^i$ into $h^i$ with a different $\tilde{\vartheta}_i$, the sufficient condition is the same with $\vartheta_i$ replaced by $\tilde{\vartheta}_i$. Our results are considering the interplay between $\xi^i$ and $\theta_i, \vartheta_i$. One can decompose the generator by choosing different functions $f^i$ and $h^i$ such that one can have a weaker condition which guarantee the existence of solution.
%\end{remark}

\begin{appendix}
\section{An auxiliary result for one-dimensional BSDEs}\label{A}

In this section, we provide an extension of \cite[Lemma 2.5]{HT}. Consider the one-dimensional BSDE
\begin{equation}\label{2}
Y_t=\xi+\int_{t}^{T}[f(s,Z_s)+h_s]ds-\int_{t}^{T}Z_sdW_s, \quad t\in[0,T],
\end{equation}
where $f:\Omega\times[0,T]\times\mathbb{R}^d\rightarrow\mathbb{R}$ is $\mathcal{P}\otimes \mathcal{B}(\mathbb{R}^d)$-measurable and $h:\Omega\times[0,T] \to \mathbb{R}$ is $\mathcal{P}$-measurable.
We assume that $f(\omega,t,z)$ is continuous in $z$ for $P\otimes dt$-almost all $(\omega,t)\in\Omega\times [0,T]$ and there exists a constant $\gamma\geq 0$ such that $$|f(\cdot,z)|\leq \gamma |z|^2, \quad \text{for all }z\in \mathbb{R}^d.$$
Consider the following conditions:
\begin{itemize}
\item[(B1)] There exists a constant $\theta\geq 0$ such that
\begin{equation*}
|f(\cdot,z)-f(\cdot,\bar{z})|\leq \theta(|z|+|\bar{z}|)|z-\bar{z}|, \quad \text{for all } z,\bar{z}\in\mathbb{R}^d.
\end{equation*}
\item[(B2)] $E[e^{2\gamma|\xi+\int_{0}^{T}h_sds|}]<\infty$.
\item[(B3)] $\xi\in L^{\infty}(\mathcal{F}_T)$ and $|h_s|\leq |z_s|^{2}$ where $z\cdot W$ is a BMO martingale with $\|z\cdot W\|_{BMO}<\frac{1}{\sqrt{2\gamma}}$.
\item[(B4)] $\xi=E[\xi]+\int_{0}^{T}v_sdW_s$ and $|h_s|\leq |z_s|^{2}$ where $v\cdot W$ and $z\cdot W$ are BMO martingales such that $\|v\cdot W\|_{BMO}<\frac{1}{16\gamma}$ and $\|z\cdot W\|_{BMO}<\frac{1}{\sqrt{4\gamma}}$.
\end{itemize}

\begin{lemma}\label{le}
\begin{itemize}[fullwidth]
\item[(i)] If (B1) holds, then the BSDE \eqref{2} has at most one solution $(Y,Z)$ such that $Z\cdot W$ is a BMO martingale.
\item[(ii)] If (B2) holds, then the BSDE \eqref{2} has a solution $(Y,Z)$ such that
 \begin{equation}\label{3}
 -\frac{1}{2\gamma}\log E\left[e^{-2\gamma(\xi+\int_{t}^{T}h_sds)}\big|\mathcal{F}_{t}\right]\leq Y_t\leq \frac{1}{2\gamma}\log E\left[e^{2\gamma(\xi+\int_{t}^{T}h_sds)}\big|\mathcal{F}_{t}\right],
 \end{equation}
 for all $t\in[0,T]$.
 \item[(iii)] If (B3) holds, then the BSDE \eqref{2} has a solution $(Y,Z)$ such that $Y\in\mathcal{S}^{\infty}(\mathbb{R})$ and $Z\cdot W$ is a BMO martingale with
\begin{align*}
&\frac{1}{2} \|Z\cdot W\|^2_{BMO} \\ \leq & \|E[u\left(\gamma,|\xi|\right)|\mathcal{F}_t]\|_{\infty}+\bigg\|u^\prime\bigg(\gamma,\frac{1}{2\gamma}\log\bigg(\frac{E[e^{2\gamma|\xi|}|\mathcal{F}_t]}{1-2\gamma\|z\cdot W\|^{2}_{BMO}}\bigg)\bigg)\bigg\|_{\infty}\|z\cdot W\|_{BMO}^2.
\end{align*}
 \item[(iv)] If (B4) holds, then the BSDE \eqref{2} has a solution $(Y,Z)$ such that $Z\cdot W$ is a BMO martingale with
 \begin{align*}
\frac{1}{2}\|Z\cdot W\|^2_{BMO} \leq  \left(\frac{1}{1-16\gamma\|v\cdot W\|_{BMO}}+\frac{1}{1-4\gamma\|z\cdot W\|^2_{BMO}}\right)\|v\cdot W\|^2_{BMO}.
\end{align*}
\end{itemize}
\end{lemma}

\begin{proof}
\begin{itemize}
\item[(i)] Let $(Y,Z)$ and $(\tilde{Y},\tilde{Z})$ be two solutions of the BSDE \eqref{2} such that $Z\cdot W$ and $\tilde{Z}\cdot W$ are BMO martingales. Denote $\Delta Y:=\tilde{Y}-Y$ and $\Delta Z:=\tilde{Z}-Z$. Then one has
\begin{align*}
\Delta Y_t&=\int_{t}^{T}(f(s,\tilde{Z}_s)-f(s,Z_s))ds-\int_{t}^{T}\Delta Z_sdW_s\\
&=\int_{t}^{T}(b_s\Delta Z_s)ds-\int_{t}^{T}\Delta Z_sdW_s=-\int_{t}^{T}\Delta Z_sd\tilde{W}_s
\end{align*}
where $\tilde{W}_t := W_t-\int_0^t b_sds$ and $b$ satisfies $|b_s|\leq\theta(|Z_s|+|\tilde{Z}_s|)$ which defines an equivalent probability measure $\tilde{P}$ by $\frac{d\tilde{P}}{dP}:=\mathcal{E}_T( b\cdot W)$. By taking the conditional expectation with respect to $\tilde{P}$ and $\mathcal{F}_t$, one has $\Delta Y=0$.
\item[(ii)] By \cite[Theorem 2]{BH}, the BSDE
\begin{equation*}
\bar{Y}_t=\left(\xi+\int_{0}^{T}h_sds \right)+\int_{t}^{T}f(s,\bar{Z}_s)ds-\int_{t}^{T}\bar{Z}_sdW_s, \quad t\in[0,T],
\end{equation*}
admits a solution $(\bar{Y},\bar{Z})$ such that
 \begin{equation*}
 -\frac{1}{2\gamma}\log E\left[e^{-2\gamma(\xi+\int_{0}^{T}h_sds)}\big|\mathcal{F}_{t}\right]\leq \bar{Y}_t\leq \frac{1}{2\gamma}\log E\left[e^{2\gamma(\xi+\int_{0}^{T}h_sds)}\big|\mathcal{F}_{t}\right].
 \end{equation*}
 Defining $Y_t:=\bar{Y}_t - \int_0^t h_s ds$, the pairing $(Y,\bar{Z})$ satisfies the BSDE \eqref{2} and $Y$ satisfies \eqref{3}.
 \item[(iii)] By the John-Nirenberg inequality \cite[Theorem 2.2]{Ka}, one has $E[\exp(2\gamma\int_{t}^{T}z_{s}^2ds)|\mathcal{F}_t]<\infty$ for all $t\in [0,T]$. Therefore it follows from \eqref{3} that
\begin{equation}\label{bddY}
|Y_t|\leq\frac{1}{2\gamma}\log\left(\frac{E[e^{2\gamma|\xi|}|\mathcal{F}_t]}{1-2\gamma\|z\cdot W\|^{2}_{BMO}}\right).
\end{equation}
By It\^{o}'s formula and the growth condition on $f$ and (A3) it holds that
\begin{align*}
u(\gamma,|Y_{t}|)&= u(\gamma,|\xi|)-\int_{t}^{T}u^{\prime}(\gamma,|Y_{s}|)\text{sgn}(Y_{s})Z_{s}dW_s\\
&+\int_{t}^{T}\bigg(u^\prime(\gamma,|Y_{s}|)\text{sgn}(Y_s)(f(s,Z_{s})+h_s)-\frac{1}{2}u^{\prime\prime}(\gamma,|Y_{s}|)|Z_s|^2\bigg)ds\\
&\leq u(\gamma,|Y_{T}|)-\int_{t}^{T}u^{\prime}(\gamma,|Y_{s}|)\text{sgn}(Y_{s})Z_{s}dW_s\\ &+\int_{t}^{T}u^\prime(\gamma,|Y_{s}|)|z_{s}|^2ds -\frac{1}{2}\int_{t}^{T}|Z_s|^2ds.
\end{align*}
Taking conditional expectation with respect to $\mathcal{F}_t$, using \eqref{bddY} and (B3), one obtains
\begin{align*}
\frac{1}{2} E\left[\int_{t}^{T}|Z_s|^2ds \big|\mathcal{F}_t\right] &\leq E[u\left(\gamma,|\xi|\right)|\mathcal{F}_t]\\ &+ E\left[\int_t^Tu^\prime\left(\gamma,\frac{1}{2\gamma}\log\left(\frac{E[e^{2\gamma|\xi|}|\mathcal{F}_s]}{1-2\gamma\|z\cdot W\|^{2}_{BMO}}\right)\right)z^2_sds\bigg|\mathcal{F}_t\right].
\end{align*}
\item[(iv)] By the John-Nirenberg inequality \cite[Theorem 2.1]{Ka}, it follows from $\|v\cdot W\|_{BMO}<\frac{1}{16\gamma}$ that
\begin{equation*}
E[e^{4\gamma|\xi|}]\leq \frac{e^{4\gamma |E[\xi]|}}{1-16\gamma\|v\cdot W\|_{BMO}}.
\end{equation*}
Combining the previous estimate with \eqref{3} we conclude that $\hat{Y}_t:=Y_t-E[\xi|\mathcal{F}_t]\in \mathcal{S}^\infty(\mathbb{R})$. Indeed we have
\begin{equation*}
|\hat{Y}_t|\leq\frac{1}{4\gamma}\log\left(\frac{1}{1-16\gamma\|v\cdot W\|_{BMO}}\right)+\frac{1}{4\gamma}\log\left(\frac{1}{1-4\gamma\|z\cdot W\|^2_{BMO}}\right).
\end{equation*}
Moreover, $\hat{Y}$ satisfies the BSDE
\begin{equation*}
\hat{Y}_t=\int_t^T (f(s,Z_s)+h_s)ds- \int_t^T (Z_s-v_s)dW_s, \quad t\in[0,T].
\end{equation*}
Applying It\^{o}'s formula and arguing as in (iii), and using additionally the inequality $(a-b)^2\geq \frac{1}{2}b^2-a^2$, one obtains for all $0\leq t\leq T$
\begin{align*}
u(2\gamma,|\hat{Y}_{t}|)&=-\int_{t}^{T}u^{\prime}(2\gamma,|\hat{Y}_{s}|)\text{sgn}(\hat{Y}_{s})(Z_{s}-v_s)dW_s \\ &+\int_{t}^{T}\bigg(u^\prime(2\gamma,|\hat{Y}_{s}|)\text{sgn}(\hat{Y}_s)(f(s,Z_{s})+h_s)-\frac{1}{2}u^{\prime\prime}(2\gamma,|\hat{Y}_{s}|)|Z_s-v_s|^2\bigg)ds\\
&\leq -\int_{t}^{T}u^{\prime}(2\gamma,|\hat{Y}_{s}|)\text{sgn}(\hat{Y}_{s})(Z_{s}-v_s)dW_s \\ &+\int_{t}^{T}\bigg(u^\prime(2\gamma,|\hat{Y}_{s}|)\text{sgn}(\hat{Y}_s)(f(s,Z_{s})+h_s)+\frac{1}{2}u^{\prime\prime}(2\gamma,|\hat{Y}_{s}|)|v_s|^2-\frac{1}{4}u^{\prime\prime}(2\gamma,|\hat{Y}_{s}|)|Z_s|^2\bigg)ds.
\end{align*}
By a similar argument as in (iii), it holds that
\begin{align*}
&\frac{1}{2}\|Z\cdot W\|^2_{BMO}\\
&\leq u''\left(2\gamma,\frac{1}{4\gamma}\log\left(\frac{1}{1-16\gamma\|v\cdot W\|_{BMO}}\right)+\frac{1}{4\gamma}\log\left(\frac{1}{1-4\gamma\|z\cdot W\|^2_{BMO}}\right)\right)\|v\cdot W\|^2_{BMO}\\
&\leq\left(\frac{1}{1-16\gamma\|v\cdot W\|_{BMO}}+\frac{1}{1-4\gamma\|z\cdot W\|^2_{BMO}}\right)\|v\cdot W\|^2_{BMO}.
\end{align*}
\end{itemize}
\end{proof}
\end{appendix}


\begin{thebibliography}{99}                                                                                               %
\bibitem{B} J.M. Bismut. Conjugate convex functions in optimal stochastic control. \textit{J. Math. Anal. Appl.} 44: 384-404, 1973.
\bibitem{BEH} K. Bahlali, E. H. Essaky, M. Hassani. Multidimensional BSDEs with super-linear growth coefficient: Application to degenerate systems of semilinear PDEs. \textit{C.R. Acad. Sci. Paris S\'er. I Math.} 348(11-12): 677-682, 2010.
%\bibitem{BEO} K. Bahlali, M. Eddahbib, Y. Oukninec. Solvability of some quadratic BSDEs without exponential moments. \textit{C.R. Acad. Sci. Paris S\'er. I Math.} 351(5-6): 229-233, 2013.
%\bibitem{BE} P. Briand and R. Elie. A simple constructive approach to quadratic BSDEs with or without delay. \textit{Stochastic Process. Appl.} 123: 2921-2939, 2013.
\bibitem{BH} P. Briand and Y. Hu. BSDE with quadratic growth and unbounded terminal value. \textit{Probability Theory and Related Fields.} 136: 604-618, 2006.
\bibitem{BH1} P. Briand and Y. Hu. Quadratic BSDEs with convex generators and unbounded terminal conditions. \textit{Probability Theory and Related Fields.} 141: 543-567, 2008.
\bibitem{CN1} P. Cheridito and K. Nam. BSDEs with terminal conditions that have bounded Malliavin derivative. \textit{Journal of Functional Analysis} 266: 1257-1285, 2014.
\bibitem{CN} P. Cheridito and K. Nam. Multidimensional quadratic and subquadratic BSDEs with special structure. \textit{Stochastics} 87: 871-884, 2015.
\bibitem{CM} B. Chikvinidze and M. Mania. New proofs of some results on bounded mean oscillation martingales
using Backward stochastic differential equations. \textit{J. Theor. Probab.} 27: 1213-1228, 2014.
\bibitem{DHB} F. Delbaen, Y. Hu, X. Bao. Backward SDEs with superquadratic growth. \textit{Probability Theory and Related Fields.} 150: 145-192, 2011.
\bibitem{EH} N. El Karoui and S. Hamad\`{e}ne. BSDEs and risk-sensitive control, zero-sum and non-zero sum game problems of stochastic functional differential equations. \textit{Stochastic Process. Appl.} 107: 145-169, 2003.
\bibitem{F} C. Frei. Splitting multidimensional BSDEs and finding local equilibria. \textit{Stochastic Process. Appl.} 124: 2654-2671, 2014.
\bibitem{FR} C. Frei and dos Reis. A financial market with interacting investors: does an equilibrium exist? \textit{Math. Financ. Econ.} 4: 161-182, 2011.
\bibitem{HT} Y. Hu and S. Tang. Multi-dimensional backward stochastic differential equations of diagonally quadratic generators. \textit{Stochastic Process. Appl.} 126: 1066-1086, 2016.
\bibitem{HT1} Y. Hu and S. Tang. Non-zero sum quadratic differential game of BSDEs and multi-dimensional diagonally quadratic BSDE. \textit{IFAC-PapersOnline} 49: 308-309, 2016.
\bibitem{KXZ} C. Kardaras, H. Xing and G. \v{Z}itkovi\'{c}. Incomplete stochastic equilibria with exponential utilities close to Pareto optimality. \textit{ArXiv-Preprint.} 1505.07224v1, 2015.
\bibitem{Ka} N. Kazamaki. Continuous Exponential Martingale and BMO. Volume 1579 of \textit{Lecture Notes in Mathematics}. Springer-Verlag, Berlin, 1994.
\bibitem{Ko} M. Kobylanski. Backward stochastic differential equations and partial differential equations with quadratic growth. \textit{Annals of Probability.} 28(2): 558-602, 2000.
\bibitem{KP} D. Kramkov and S. Pulido. A system of quadratic BSDEs arising in a price impact model. \textit{Ann. Appl. Probab.} 26: 794-817, 2016.
\bibitem{MR} F. Masiero and A. Richou. A note on the existence of solutions to Markovian superquadratic BSDEs with an unbounded ternimal condtion. \textit{Electron. J. Probab} 18(50): 1-15, 2013.
\bibitem{PP} E.Pardoux and S. G. Peng. Adapted solution of a backward stochastic differential equations. \textit{System Controll Lett.} 14: 55-61, 1990.
\bibitem{R} A. Richou. Markovian quadratic and superquadratic BSDEs with an unbounded terminal condition. \textit{Stochastic Process. Appl.} 122: 3173-3208, 2012.
\bibitem{T} R. Tevzadze. Solvability of backward stochastic differential equations with quadratic growth. \textit{Stochastic Process. Appl.} 118: 503-515, 2008.
\bibitem{NT} N. Touzi. Optimal Stochastic Control, Stochastic Target Problems, and Backward SDE. In: Fields Institute Monographs, vol. 29. Springer, New York, 2013.
\bibitem{XZ} H. Xing and G. \v{Z}itkovi\'{c}. A class of globally solvable Markovian quadratic BSDE systems and applications. \textit{ArXiv-Preprint.} 1603.00217, 2016.

\end{thebibliography}
\end{document}